\documentclass{aptpub}

\authornames{D. Hong, S. Man, J-C. Birget, D. S. Lun} 
\shorttitle{A wavelet-based approximation of fBm} 

\begin{document}

\title{A wavelet-based almost sure uniform approximation of fractional Brownian motion with a parallel algorithm}

\authorone[Rutgers University]{Dawei Hong}
\addressone{Center for Computational and Integrative Biology, Department of Computer Science, Rutgers University, Camden, NJ 08102, USA; Email address: dhong@camden.rutgers.edu}
\authortwo[Southwest Minnesota State University]{Shushuang Man}
\addresstwo{Department of Mathematics and Computer Science, Southwest Minnesota State University, Marshall, MN 56258, USA}
\authorone[Rutgers University]{Jean-Camille Birget}
\authorone[Rutgers University]{Desmond S. Lun}

\begin{abstract}
We construct a wavelet-based almost sure uniform approximation of fractional Brownian motion (fBm) $(B_t^{(H)})_{t\in[0, 1]}$ of Hurst index $H$ $\in$ $(0, 1)$. Our results show that by Haar wavelets which merely have one vanishing moment, an almost sure uniform expansion of fBm of $H$ $\in$ $(0, 1)$ can be established. The convergence rate of our approximation is derived. We also describe a parallel algorithm that generates sample paths of an fBm  efficiently.
\end{abstract}

\smallskip



\section{Introduction}
\label{Intro}
A fractional Brownian motion (fBm) $(B_t^{(H)})_{t\in[0, T]}$ of Hurst index $H \in (0, 1)$ is a centered Gaussian process with covariance $E[B_{t_1}^{(H)} B_{t_2}^{(H)}]$ $=$ $(1/2)(t_1^{2H}+t_2^{2H}-|t_1-t_2|^{2H})$ for all $t_1, t_2\in[0, T]$. A standard Brownian motion (Bm) $(B_t)_{t\in[0, T]}$ is the special case of $H$ $=$ $1/2$. There are a great number of applications of fBm in engineering and the sciences; see \cite{BHOZ} and references therein. The study of approximations of fBm has been active since the 1970s. A major focus is to find approximations of fBm that converge in law; for example, see \cite{Dav, Ta, DJ, BJT, LD} and references therein. However, practical implementations often require almost sure uniform, also termed as strong uniform, approximations of fBm. It works as follows: Let $(B_t^{(H)})_{t\in[0, 1]}$ be an fBm of some $H$ $\in$ $(0, 1)$. Then, with respect to the probability space where $(B_t^{(H)})_{t\in[0, 1]}$ is defined, the following event occurs with probability $1$. For a sample path of $(B_t^{(H)})_{t\in[0, 1]}$ there is a sequence of functions of $t$ $\in$ $[0, 1]$ produced by the approximation which uniformly converges to the sample path; conversely, a sequence of functions of $t$ $\in$ $[0, 1]$ produced by the approximation uniformly converges to a sample path of $(B_t^{(H)})_{t\in[0, 1]}$.

Meyer, Sellan and Taqqu \cite{MST} obtained several wavelet series expansions of fBm of $H$ $\in$ $(0, 1)$ that almost surely and uniformly converge. Their results brought deep insights into spectral properties of fBm. For instance, the wavelet series expansion of fBm in section 7 yielded a very, if not the most, efficient mathematical representation of the spectral properties of fBm -- a subject that has attracted much research for decades -- showed in section 8.

K\"uhn and Linde \cite{KL} showed that the optimal convergence rate that a series expansion of fBm may reach is $O(N^{-H} \sqrt{\log N})$ if the expansion converges almost surely and uniformly. Ayache and Taqqu \cite{AT} proved that under certain conditions wavelet series expansions of fBm in \cite{MST} converge at the optimal rate. Dzhaparidze and van Zanten \cite{DvZ} constructed a series expansion of fBm of $H$ $\in$ $(0, 1)$ (in the frequency domain) which almost surely and uniformly converges at the optimal rate \cite{DvZ1}.

The above results will have a long-lasting impact on the study of fBm; in the meantime they stimulate further studies. Theorem 2 in \cite{MST} and remark 4 on the theorem motivated our investigation. Haar wavelets are very convenient to compute. Moreover, the simple form of the Mandelbrot - van Ness representation of fBm \cite{MV} is likely to yield a fast algorithm. A question is whether using the Mandelbrot - van Ness representation and Haar wavelets, one can construct an almost sure uniform approximation of fBm for all $H$ $\in$ $(0, 1)$. In this paper, we establish such an approximation of fBm of $H$ $\in$ $(0, 1)$. Our approach is to apply L\'evy's equivalence theorem (e.g., see Theorem 9.7.1 in \cite{Du}) to a Haar wavelet-based approximation of fBm obtained from the Mandelbrot - van Ness representation, and then to carefully evaluate the wavelet coefficients.

As shown by \cite{MST}, wavelet approximation of fBm is a powerful approach. A key idea of this approach is to almost surely and uniformly approximate the sample paths in a process, using i.i.d. Gaussian random variables with a finely designed basis of $L^2$ space such as Meyer's or Daubechies' wavelets. The conditions for wavelet approximations of fBm with the optimal convergence rate \cite{AT} need wavelets to have the first six vanishing moments. It is a question if one can use Haar wavelets that merely have the first vanishing moment to obtain an almost sure uniform approximation of fBm for all $H$ $\in$ $(0, 1)$. We show this is possible. The convergence rate of our almost sure uniform approximation of fBm by Haar wavelets reaches the optimal $O(N^{- H} \sqrt{\log N})$ for $H$ $\in$ $(0, 1/2]$, but the convergence slows down to rate $O(N^{-(1 - H)} \sqrt{\log N})$ for $H$ $\in$ $(1/2, 1)$ (Theorem \ref{Thms2}). Haar wavelets (piecewise constant functions) do not introduce computational errors by themselves; and our approximation (based on the Mandelbrot - van Ness representation) is in a rather simple form. These two advantages make our approximation of fBm suitable for practical applications when $H$ is not close to $1$. We also describe a parallel algorithm that generates sample paths of an fBm efficiently.

Section \ref{Sec2} is for preliminaries. In sections \ref{Sec3}, \ref{Sec4}, \ref{Sec5} and \ref{Sec6}, we construct and prove an almost sure uniform approximation of fBm of $H$ $\in$ $(0, 1)$. We describe a parallel algorithm for approximation of fBm in section \ref{Sec7}.

\section{Preliminaries}
\label{Sec2}
Let $C_H$ $=$ $(\Gamma(H + 1/2))^{-1}$, the reciprocal of the Gamma function at $H + 1/2$. The Mandelbrot - van Ness stochastic integral representation of fBm \cite{MV} is $B_t^{(H)}$ $=$ $C_H\int_{-\infty}^t \left( (t-s)_+^{H-1/2}-(-s)_+^{H-1/2} \right)dB_s$ for $H$ $\in$ $(0, 1/2)$ $\cup$ $(1/2, 1)$; and when $H$ $=$ $1/2$, fBm becomes Bm. In what follows, we denote the underlying probability space for the above representation of fBm by $(\Omega, {\cal F}, \Pr)$ where $\cal F$ is a standard Brownian filtration. Our construction of an almost sure uniform approximation of fBm is based on a rewriting of the Mandelbrot - van Ness stochastic integral representation
\begin{equation}
\label{Sec2Eq1}
\begin{split}
B_t^{(H)} = I_1(t, H) +  I_2(t, H) + I_3(t, H), \ t \in [0, 1], \ \text{ where}
\end{split}
\end{equation}
$I_1(t, H)$ $=$ $C_H\int_{0}^t(t-s)^{H-1/2}dB_s$, $I_2(t, H)$ $=$ $C_H\int_{-1}^0 \left((t-s)^{H-1/2}-(-s)^{H-1/2}\right)dB_s$, and $I_3(t, H)$ $=$ $C_H\int_{-\infty}^{-1} \left((t-s)^{H-1/2}-(-s)^{H-1/2}\right)dB_s$.

Let $(\phi_n)_{n \geq 0}$ be a complete orthonormal basis for $L^2[a, b]$. For $f$ $\in$ $L^2[a, b]$, we have $f$ $=$ $\sum_{n = 0}^{\infty} \langle f, \phi_n \rangle \phi_n$ in $L^2[a, b]$. On both sides of $f$ $=$ $\sum_{n = 0}^{\infty} \langle f, \phi_n \rangle \phi_n$, we take Wiener integration. Then we informally interchange the order of integration and summation on the right side, having $\int_a^b f(s) dB_s$ $=$ $\sum_{n = 0}^{\infty} \langle f, \phi_n \rangle \int_a^b \phi_n(s) dB_s$. By L\'evy's equivalence theorem we have
\begin{thm}
\label{Sec2Thms1}
$\lim_{N \rightarrow \infty} \sum_{n = 0}^N \langle f, \phi_n \rangle \int_a^b \phi_n(s) dB_s$ $=$ $\int_a^b f(s) dB_s$ almost surely. $\Box$
\end{thm}

The Haar wavelet on $[0, 1]$ is defined as follows. Let ${\cal H}(s)$ $=$ $1$ if $s$ $\in$ $[0, 1/2)$, ${\cal H}(s)$ $=$ $-1$ if $s$ $\in$ $[1/2, 1]$, and ${\cal H}(s)$ $=$ $0$ otherwise. For $n$ $=$ $2^j$ $+$ $k$ with $j$ $\geq$ $0$ and $0 \leq k < 2^j$, define ${\cal H}_n(s)$ $=$ $2^{j/2}{\cal H}(2^js-k)$ and ${\cal H}_0(s) = 1$. The sequence $({\cal H}_n)_{n\geq0}$ is the Haar wavelet on $[0,1]$, which constitutes a complete orthonormal basis for $L^2[0, 1]$. In a similar way, we can define the Haar wavelet on any given interval $[a, b]$ $\subset$ $\mathbb R$ to constitute a complete orthonormal basis for $L^2[a, b]$ (see \cite{Da}).

\section{Approximation of $I_1(t, H)$}
\label{Sec3}
We construct and prove an almost sure uniform approximation of $I_1(t, H)$. Consider a family of functions $f_t^{(1)}$ $\in$ $L^2[0,1]$ with a parameter $t$ $\in$ $(0, 1] \cap \mathbb Q$:
\begin{equation*}
\begin{split}
&f_t^{(1)}(s) =
\begin{cases}
(t-s)^{H-1/2} & {\rm if} \ s \in [0, t) \\
0 & \text{otherwise.}
\end{cases}\\
\end{split}
\end{equation*}
By Theorem \ref{Sec2Thms1} we have
\begin{equation}
\label{Sec3Eq1}
\begin{split}
\Pr\left\{ \left(\int_0^1 f_t^{(1)}(s) dB_s\right)(\omega) = \left(\sum_{n = 0}^{\infty} \langle f_t^{(1)}, {\cal H}_n \rangle \int_0^1 {\cal H}_n (s) dB_s\right)(\omega) \right\} = 1
\end{split}
\end{equation}
for each $t$ $\in$ $(0, 1] \cap \mathbb Q$ and as a consequence,
\begin{equation*}
\begin{split}
\Pr \bigcap_{t \in (0, 1] \cap \mathbb Q} & \left\{ \left(\int_0^1
f_t^{(1)}(s) dB_s\right)(\omega) = \left(\sum_{n = 0}^{\infty} \langle
f_t^{(1)}, {\cal H}_n  \rangle \int_0^1 {\cal H}_n(s) dB_s\right)(\omega)
\right\} = 1 .
\end{split}
\end{equation*}
We define for all $N \geq 1$,
\begin{equation}
\label{Sec3Eq2}
\begin{split}
W_1(t, H, N) =
\begin{cases}
C_H \sum_{n=0}^N \langle f_t^{(1)}, {\cal H}_n \rangle {\cal L}_n^{(1)}
& \text{ for } t\in (0, 1] \cap \mathbb Q\\
0 & \text{ for } t = 0.
\end{cases}
\end{split}
\end{equation}
Here ${\cal L}_n^{(1)}$ $=$ $\int_0^1 {\cal H}_n(s) dB_s$, $n$ $=$ $0$, $1$, $\ldots$, $N$, are i.i.d. Gaussian random variables with mean $0$ and variance $1$.

In what follows, we use two conventions: $n \in \mathbb Z^+$ is said to be at level $j$ if $n$ $=$ $2^j$ $+$ $k$ with $j$ $\geq$ $0$ and $0$ $\leq$ $k$ $<$ $2^j$; and the interval $[\frac{k}{2^j}, \frac{k+1}{2^j})$ is meant to be $[\frac{k}{2^j}, \frac{k+1}{2^j}]$ when $\frac{k+1}{2^j}$ $=$ $1$.
\begin{lem}
\label{Sec3Lems1}
There is an absolute constant $D_1$ $>$ $0$ such that for every $t$ $\in$ $(0, 1] \cap \mathbb Q$ and for all $N$ $>$ $1$,
$\sum_{n = N + 1}^{\infty}\left\langle f_t^{(1)}, {\cal H}_n \right\rangle^2$ $\leq$ $D_1 \left(H(1-H)N^{2H}\right)^{-1}$.
\end{lem}
\textbf{Proof.}
For $t$ $\in$ $(0, 1] \cap \mathbb Q$, at each level $j$ $=$ $0$, $1$, $\ldots$, we partition the set
\begin{equation*}
\begin{split}
\left\{n = 2^j + k : k = 0,1,\ldots,2^j - 1\right\}
\end{split}
\end{equation*}
into three subsets: ${\cal G}_1(j, t)$ consisting of all $n$ ($=$ $2^j + k$) such that $[\frac{k}{2^j}, \frac{k+1}{2^j})$ $\subseteq$ $[0, t)$; ${\cal G}_2(j, t)$ consisting of the one $n$ such that $t$ $\in$ $[\frac{k}{2^j}, \frac{k+1}{2^j})$; and ${\cal G}_3(j, t)$ consisting of all $n$ such that $t$ $\notin$ $\bigcup_{k^* = k}^{2^j - 1} [\frac{k^*}{2^j}, \frac{k^* + 1}{2^j})$.

Consider a fixed $j$. By the definition of $f_t^{(1)}$ we have
\begin{equation}
\label{lems2-1-1-01Eq1}
\begin{split}
\left\langle f_t^{(1)}, {\cal H}_n \right\rangle = 0 \ \text{ for every } n \in {\cal G}_3(j, t).
\end{split}
\end{equation}

For the only $n$ $\in$ ${\cal G}_2(j, t)$, we denote by $\widehat{k_{t, j}}$ the $k$ that appears in $n$ $=$ $2^j + k$. We have
\begin{equation*}
\begin{split}
&\left\langle f_t^{(1)}, {\cal H}_n \right\rangle = 2^{j/2} \left[\int_{\frac{2\widehat{k_{t, j}}}{2^{j + 1}}}^{\frac{2\widehat{k_{t, j}} + 1}{2^{j + 1}}} f_t^{(1)}(s) ds - \int_{\frac{2\widehat{k_{t, j}} + 1}{2^{j + 1}}}^{\frac{2\widehat{k_{t, j}} + 2}{2^{j + 1}}} f_t^{(1)}(s) ds \right]
\end{split}
\end{equation*}
which implies
\begin{equation*}
\begin{split}
&\left|\left\langle f_t^{(1)}, {\cal H}_n \right\rangle\right| \leq 2^{j/2} \times \\
&\max \left\{\int_{\frac{2\widehat{k_{t, j}}}{2^{j + 1}}}^{\frac{2\widehat{k_{t, j}} + 1}{2^{j + 1}}} \left(\frac{2\widehat{k_{t, j}} + 1}{2^{j + 1}} - s \right)^{H - 1/2} ds, \int_{\frac{2\widehat{k_{t, j}} + 1}{2^{j + 1}}}^{\frac{2\widehat{k_{t, j}} + 2}{2^{j + 1}}} \left(\frac{2\widehat{k_{t, j}} + 2}{2^{j + 1}} - s \right)^{H - 1/2} ds \right\}.
\end{split}
\end{equation*}
Using this inequality, by calculation we have for $n$ $\in$ ${\cal G}_2(j, t)$,
\begin{equation}
\label{lems2-1-1-01Eq2}
\begin{split}
\left\langle f_t^{(1)}, {\cal H}_n \right\rangle^2 \leq 2^{-2jH} \left(2^{-(2H + 1)} \, (H + 1/2)^{-2}\right).
\end{split}
\end{equation}

For each $n$ ($=$ $2^j + k$) $\in$ ${\cal G}_1(j,t)$, we have
\begin{equation}
\label{lems2-1-1-01Eq3}
\begin{split}
&\left\langle f_t^{(1)}, {\cal H}_n \right\rangle = 2^{j/2} \left[\int_{\frac{2k}{2^{j + 1}}}^{\frac{2k + 1}{2^{j + 1}}} (t - s)^{H - 1/2} ds - \int_{\frac{2k + 1}{2^{j + 1}}}^{\frac{2k + 2}{2^{j + 1}}} (t - s)^{H - 1/2} ds \right]\\
& = \frac{2^{j/2}}{H + 1/2}\left[\left(\left(t-\frac{2k}{2^{j+1}}\right)^{H+1/2} - \left(t-\frac{2k+1}{2^{j+1}}\right)^{H+1/2}\right) - \right.\\
&\left.\left(\left(t-\frac{2k + 1}{2^{j+1}}\right)^{H+1/2} - \left(t-\frac{2k+2}{2^{j+1}}\right)^{H+1/2}\right)\right].
\end{split}
\end{equation}
To facilitate our argument, we introduce a function $w$ of $h$: $w(h)$ $=$ $g(x_0 + h)$ $+$ $g(x_0 - h)$ $-$ $2g(x_0)$ where $g(\cdot)$ $=$ $\left(\cdot\right)^{H+1/2}~\text{ and } x_0$ $=$ $t-\frac{2k+1}{2^{j+1}}$. We let $h$ $=$ $\frac{1}{2^{j + 1}}$ and rewrite \eqref{lems2-1-1-01Eq3} as
\begin{equation}
\label{lems2-1-1-01Eq4}
\begin{split}
&\left\langle f_t^{(1)}, {\cal H}_n \right\rangle = \frac{2^{j/2}}{H + 1/2} \, w(h).
\end{split}
\end{equation}
By Taylor's expansion,
\begin{equation*}
\begin{split}
&w(h) = w(0) + \frac{w'(0)}{1!}h + \frac{w''(\theta h)}{2!}h^2~\text{ (for some } 0 < \theta < 1)\\
& = \frac{w''(\theta h)}{2!}h^2~\text{(since } w(0) = w'(0) = 0);
\end{split}
\end{equation*}
hence we have
\begin{equation*}
\begin{split}
&w(h) = h^2 \, \frac{w''(\theta h)}{2!} = 2^{-2(j + 1)}\frac{(H + 1/2)(H - 1/2)}{2} \times \\
&\left[\left(t-\frac{2k + 1 + \theta}{2^{j+1}}\right)^{H - 3/2} + \left(t-\frac{2k + 1 - \theta}{2^{j+1}}\right)^{H - 3/2}\right].
\end{split}
\end{equation*}
This equality leads us to consider the case where $n$ ($=$ $2^j$ $+$ $k$) $\in$ ${\cal G}_1(j, t)$ with $k$ $+$ $2$ $\leq$ $\widehat{k_{t, j}}$. In this case, by \eqref{lems2-1-1-01Eq4} we have
\begin{equation*}
\begin{split}
&\left|\left\langle f_t^{(1)}, {\cal H}_n \right\rangle\right| \leq 2^{j/2}2^{-2(j + 1)}|H - 1/2| \left(t-\frac{2k + 1 + \theta}{2^{j+1}}\right)^{H - 3/2}
\end{split}
\end{equation*}
(since $0$ $<$ $\theta$ $<$ $1$ and $0$ $<$ $H$ $<$ $1$), which yields
\begin{equation*}
\begin{split}
&\left|\left\langle f_t^{(1)}, {\cal H}_n \right\rangle\right| \leq 2^{j/2}2^{-2(j + 1)}|H - 1/2| \left(\frac{2\widehat{k_{t, j}}}{2^{j + 1}} - \frac{2k + 2}{2^{j+1}}\right)^{H - 3/2}\\
&= \frac{|H - 1/2|}{4}2^{-jH}\left(\widehat{k_{t, j}} - (k + 1)\right)^{H - 3/2}.
\end{split}
\end{equation*}
Thus, for $n$ ($=$ $2^j$ $+$ $k$) $\in$ ${\cal G}_1(j, t)$ with $k$ $+$ $2$ $\leq$ $\widehat{k_{t, j}}$, we have
\begin{equation}
\label{lems2-1-1-01Eq5}
\begin{split}
&\left|\left\langle f_t^{(1)}, {\cal H}_n \right\rangle\right|^2 \leq 2^{-2jH}\left(\widehat{k_{t, j}} - (k + 1)\right)^{2H - 3}\frac{|H - 1/2|^2}{16}.
\end{split}
\end{equation}

There is one and only one $\left\langle f_t^{(1)}, {\cal H}_n \right\rangle$ with $n$ $\in$ ${\cal G}_1(j, t)$ which is not included in \eqref{lems2-1-1-01Eq5}, namely $n$ $=$ $2^j$ $+$ $\widehat{k_{t, j}} - 1$. However, in this case we have
\begin{equation*}
\begin{split}
&\left\langle f_t^{(1)}, {\cal H}_n \right\rangle = 2^{j/2} \left[\int_{\frac{2\widehat{k_{t, j}} - 2}{2^{j + 1}}}^{\frac{2\widehat{k_{t, j}} - 1}{2^{j + 1}}} (t - s)^{H - 1/2} ds - \int_{\frac{2\widehat{k_{t, j}} - 1}{2^{j + 1}}}^{\frac{2\widehat{k_{t, j}}}{2^{j + 1}}} (t - s)^{H - 1/2} ds \right]
\end{split}
\end{equation*}
and hence
\begin{equation}
\label{lems2-1-1-01Eq6}
\begin{split}
&\left\langle f_t^{(1)}, {\cal H}_n \right\rangle^2 \leq \frac{2^j}{(H + 1/2)^2}2^{-(2H + 1)j} = 2^{-2jH}(H + 1/2)^{-2}.
\end{split}
\end{equation}

Now, putting \eqref{lems2-1-1-01Eq1}, \eqref{lems2-1-1-01Eq2}, \eqref{lems2-1-1-01Eq5} and \eqref{lems2-1-1-01Eq6} together, we have that there is an absolute constant $D_1^*$ $>$ $0$ such that at any level $j$,
\begin{equation*}
\begin{split}
&\sum_{\{n \text{ at level } j\}} \left|\left\langle f_t^{(1)}, {\cal H}_n \right\rangle\right|^2 \leq D_1^* \,2^{-2jH}
\sum_{\ell = 1}^{\infty}\left(\frac{1}{\ell}\right)^{3 - 2H}\\
&= D_1^* \,2^{-2jH} \left(1 + \sum_{\ell = 2}^{\infty}\left(\frac{1}{\ell}\right)^{3 - 2H} \right) \leq D_1^* \,2^{-2jH} \left(1 + \int_1^\infty \frac{dv}{v^{3 - 2H}}\right).
\end{split}
\end{equation*}
This inequality can be written as
\begin{equation*}
\begin{split}
&\sum_{\{n \text{ at level } j\}} \left|\left\langle f_t^{(1)}, {\cal H}_n \right\rangle\right|^2 \leq \frac{D_1^{**}}{1 - H}\,2^{-2jH}, \ D_1^{**} > 0 \text{ is an absolute constant}.
\end{split}
\end{equation*}
Therefore we have
\begin{equation*}
\begin{split}
&\sum_{n = N + 1}^{\infty}\left\langle f_t^{(1)}, {\cal H}_n \right\rangle^2 \leq \sum_{j = \lfloor \log_2 N \rfloor }^{\infty}
\sum_{\{n \text{ at level } j\}} \left|\left\langle f_t^{(1)}, {\cal H}_n \right\rangle\right|^2 \leq \sum_{j = \lfloor \log_2 N \rfloor }^{\infty} \frac{D_1^{**}}{1 - H}\,2^{-2jH} \\
&= \frac{D_1^{**}}{1 - H}\,2^{-2 \lfloor \log_2 N \rfloor H} \sum_{j = 0}^{\infty} 2^{-2jH} = \frac{D_1^{**}}{1 - H} \, 2^{-2 \lfloor \log_2 N \rfloor H} \frac{1}{1 - 2^{-2H}}.
\end{split}
\end{equation*}
The lemma follows from this inequality and the fact that there is an absolute constant $G$ $>$ $0$ such that $1/(1 - 2^{-2H})$ $\leq$ $G/H$ for all $H$ $\in$ $(0, 1)$ (because $\lim_{H \rightarrow 0_+} (1 - 2^{-2H})/H$ $=$ $2\log 2$). $\Box$

\begin{lem}
\label{Sec3Lems2}
For any given $H$ $\in$ $(0,1)$ and $q \geq 2$, we have for all $N$ $>$ $1$,
\begin{equation*}
\begin{split}
&\Pr\left\{\sup_{t \in [0,1] \cap \mathbb Q}\left| I_1(t, H) - W_1(t, H, N) \right| \geq \frac{C_H\sqrt{2 D_1 q}}{\sqrt{H(1 - H)}} \frac{\sqrt{\log N}}{N^H}\right\} \leq \frac{1}{\sqrt{\pi} N^q},
\end{split}
\end{equation*}
where $D_1$ is the absolute constant used in Lemma \ref{Sec3Lems1}.
\end{lem}
\textbf{Proof.}
By definition $I_1(0, H)$ $=$ $0$ $=$ $W_1(0, H, N)$. So, we focus on the case of $t$ $\in$ $(0, 1] \cap \mathbb Q$. By \eqref{Sec3Eq2} and the consequence of \eqref{Sec3Eq1}, we have
\begin{equation}
\label{lems2-1-1Eq1}
\begin{split}
\Pr \bigcap_{t \in (0, 1] \cap \mathbb Q} & \{(I_1(t, H) - W_1(t, H, N))(\omega) =  \\
& C_H\sum_{n = N + 1}^\infty \langle f_t^{(1)}, {\cal H}_n \rangle \int_0^1
{\cal H}_n (s) dB_s(\omega) \} \ = \ 1.
\end{split}
\end{equation}
Here $\sum_{n=N+1}^{\infty} \langle f_t^{(1)}, {\cal H}_n \rangle\int_0^1{\cal H}_n(s)dB_s$ is a Gaussian random variable with mean $0$ and variance $\sum_{n=N+1}^{\infty} \langle f_t^{(1)}, {\cal H}_n \rangle^2$. We denote $\sum_{n=N+1}^{\infty} \langle f_t^{(1)}, {\cal H}_n \rangle^2$ by $\sigma_1^2(t, H, N)$. For any given $H$ $\in$ $(0,1)$ and $q$ $\geq$ $2$, we have
\begin{equation*}
\begin{split}
& \Pr\left\{\left|\sum_{n=N+1}^{\infty} \int_0^1\left\langle f_t^{(1)}, {\cal H}_n \right\rangle {\cal H}_n(s) dB_s\right| \geq \frac{\sqrt{2 D_1 q \log N}}{N^H\sqrt{H(1 - H)}}\right\} = \frac{\sqrt 2}{\sigma_1(t, H, N)\sqrt \pi} \times \\
& \int_{ \frac{\sqrt{2 D_1 q \log N}}{N^H\sqrt{H(1 - H)}} }^{\infty} \exp\left(-\frac{u^2}{2\sigma^2_1(t, H, N)}\right)du = \frac{2}{\sqrt \pi}\int_{ \frac{\sqrt{2 D_1 q \log N}}{\sqrt{2}\sigma_1(t, H, N) N^H\sqrt{H(1 - H)}} }^{\infty} e^{-v^2}dv \\
& \leq \frac{2}{\sqrt \pi}\int_{ \sqrt{q \log N} }^{\infty} e^{- v^2}dv \text{ (by Lemma \ref{Sec3Lems1})} \\
& \leq \frac{1}{\sqrt \pi}\int_{ \sqrt{q \log N} }^{\infty} 2 v e^{- v^2}dv \text{ (since $\sqrt{ q \log N} > 1$ for $q \geq 2$ and $N$ $>$ $1$)}.
\end{split}
\end{equation*}
Putting this and \eqref{lems2-1-1Eq1} together, we complete a proof for the Lemma. $\Box$

\section{Approximation of $I_2(t, H)$}
\label{Sec4}
Our construction and proof for an almost sure uniform approximation of $I_2(t, H)$ are similar to that for $I_1(t, H)$ presented in the previous section. Consider the Haar wavelet $(\widetilde {\cal H}_n)_{n \geq 0}$ on $[-1, 0]$. We consider a family of functions $f_t^{(2)}$ $\in$ $L^2[-1, 0]$ with a parameter $t$ $\in$ $[0, 1] \cap \mathbb Q$:
\begin{equation*}
\begin{split}
&f_t^{(2)}(s) =
\begin{cases}
(t - s)^{H-1/2} - (-s)^{H-1/2} & {\rm if} \ s \in [-1, 0),\\
0 & \text{otherwise.}
\end{cases}\\
\end{split}
\end{equation*}

By Theorem \ref{Sec2Thms1} we have
\begin{equation}
\label{Sec4Eq1}
\begin{split}
\Pr\left\{ \left(\int_{-1}^0 f_t^{(2)}(s) dB_s\right)(\omega) = \left(\sum_{n = 0}^{\infty} \langle f_t^{(2)}, \widetilde{\cal H}_n \rangle \int_{-1}^0 \widetilde{\cal H}_n (s) dB_s\right)(\omega) \right\} = 1
\end{split}
\end{equation}
for each $t$ $\in$ $[0, 1] \cap \mathbb Q$, and as a consequence
\begin{equation*}
\begin{split}
\Pr \bigcap_{t \in [0, 1] \cap \mathbb Q} & \left\{ \left(\int_{-1}^0
f_t^{(2)}(s) dB_s\right)(\omega) = \left(\sum_{n = 0}^{\infty} \langle
f_t^{(2)}, \widetilde{\cal H}_n  \rangle \int_{-1}^0 \widetilde{\cal H}_n(s)
dB_s\right)(\omega) \right\} = 1.
\end{split}
\end{equation*}
We define for all $N \geq 1$,
\begin{equation}
\label{Sec4Eq2}
\begin{split}
W_2(t, H, N) =
\begin{cases}
C_H \sum_{n=0}^N \langle f_t^{(2)}, \widetilde{\cal H}_n \rangle {\cal L}_n^{(2)}
& \text{ for } t\in [0, 1] \cap \mathbb Q\\
0 & \text{ for } t = 0.
\end{cases}
\end{split}
\end{equation}
Here ${\cal L}_n^{(2)}$ $=$ $\int_{-1}^0 \widetilde{\cal H}_n(s) dB_s$, $n$ $=$ $0$, $1$, $\ldots$, $N$, are i.i.d. Gaussian random variables with mean $0$ and variance $1$. Notice that the sequence $({\cal L}_n^{(2)})_{n \geq 0}$ is independent of the sequence $({\cal L}_n^{(1)})_{n \geq 0}$ used in the definition of $W_1(t, H, N)$.

\begin{lem}
\label{Sec4Lems1}
There is an absolute constant $D_2$ $>$ $0$ such that for every $t$ $\in$ $[0, 1] \cap \mathbb Q$ and for all $N$ $>$ $1$,
$\sum_{n = N + 1}^{\infty}\left\langle f_t^{(2)}, \widetilde{\cal H}_n \right\rangle^2$ $\leq$ $D_2 \left( H(1-H)N^{2H} \right)^{-1}$.
\end{lem}
\textbf{Proof.} For each $t$ $\in$ $[0, 1] \cap \mathbb Q$
\begin{equation}
\label{lems2-1-1-02Eq1}
\begin{split}
&\sum_{n = N + 1}^{\infty}\left\langle f_t^{(2)}, \widetilde{\cal H}_n \right\rangle^2 \\
&\leq 2\left(\sum_{n = N + 1}^{\infty}\left\langle (t - s)^{H-1/2}, \widetilde{\cal H}_n \right\rangle^2 + \sum_{n = N + 1}^{\infty}\left\langle (-s)^{H-1/2}, \widetilde{\cal H}_n \right\rangle^2\right).
\end{split}
\end{equation}
In the right side of \eqref{lems2-1-1-02Eq1}, by changing variables we have $\langle (t - s)^{H-1/2}, \widetilde{\cal H}_n \rangle$ $=$ $\langle (t + s)^{H-1/2}, {\cal H}_n \rangle$ and $\langle (-s)^{H-1/2}, \widetilde{\cal H}_n \rangle$ $=$ $\langle s^{H-1/2}, {\cal H}_n \rangle$. Below we estimate $\sum_{n = N + 1}^{\infty}\langle (t + s)^{H-1/2}, {\cal H}_n \rangle^2$ and $\sum_{n = N + 1}^{\infty}\langle s^{H-1/2}, {\cal H}_n \rangle^2$, respectively.

For $t$ $\in$ $(0, 1] \cap \mathbb Q$, at each level $j$, we have for each $n$ $=$ $2^j$ $+$ $k$, $k$ $=$ $0$, $\ldots$, $2^j - 1$,
\begin{equation}
\label{lems2-1-1-02Eq2}
\begin{split}
&\left\langle (t + s)^{H - 1/2}, {\cal H}_n \right\rangle \\
& = 2^{j/2} \left[\int_{\frac{2k}{2^{j + 1}}}^{\frac{2k + 1}{2^{j + 1}}} (t + s)^{H - 1/2} ds - \int_{\frac{2k + 1}{2^{j + 1}}}^{\frac{2k + 2}{2^{j + 1}}} (t + s)^{H - 1/2} ds \right]\\
& = \frac{2^{j/2}}{H + 1/2}\left[\left(\left(t + \frac{2k + 1}{2^{j+1}}\right)^{H+1/2} - \left(t + \frac{2k}{2^{j+1}}\right)^{H+1/2}\right) - \right.\\
&\left.\left(\left(t + \frac{2k + 2}{2^{j+1}}\right)^{H+1/2} - \left(t + \frac{2k + 1}{2^{j+1}}\right)^{H+1/2}\right)\right].
\end{split}
\end{equation}
To facilitate our argument, we introduce a revised version of the function $w$ of $h$ used in the proof of Lemma \ref{Sec3Lems1}. Since there will be no confusion, we denote this revised version by $w$ as follows: $w(h)$ $=$ $2g(x_0)$ $-$ $g(x_0 + h)$ $-$ $g(x_0 - h)$ where $g(\cdot)$ $=$ $\left(\cdot\right)^{H+1/2}$ and $x_0$ $=$ $t$ $+$ $\frac{2k+1}{2^{j+1}}$. We let $h$ $=$ $\frac{1}{2^{j + 1}}$ and rewrite \eqref{lems2-1-1-02Eq2} as
\begin{equation*}
\begin{split}
&\left\langle (t + s)^{H - 1/2}, {\cal H}_n \right\rangle = \frac{2^{j/2}}{H + 1/2} \, w(h).
\end{split}
\end{equation*}
Then by Taylor's expansion,
\begin{equation*}
\begin{split}
&w(h) = w(0) + \frac{w'(0)}{1!}h + \frac{w''(\theta h)}{2!}h^2~\text{ (for some } 0 < \theta < 1)\\
& = \frac{w''(\theta h)}{2!}h^2~\text{(since } w(0) = w'(0) = 0);
\end{split}
\end{equation*}
hence we have
\begin{equation*}
\begin{split}
&w(h) = h^2 \, \frac{w''(\theta h)}{2!} = - 2^{-2(j + 1)}\frac{(H + 1/2)(H - 1/2)}{2} \times \\
&\left[\left(t + \frac{2k + 1 + \theta}{2^{j+1}}\right)^{H - 3/2} + \left(t + \frac{2k + 1 - \theta}{2^{j+1}}\right)^{H - 3/2}\right].
\end{split}
\end{equation*}
Putting this equality, the rewriting of \eqref{lems2-1-1-02Eq2}, and the fact that $0$ $<$ $\theta$ $<$ $1$ and $0$ $<$ $H$ $<$ $1$ together, we have
\begin{equation}
\label{lems2-1-1-02Eq3}
\begin{split}
&\left|\left\langle (t + s)^{H - 1/2}, {\cal H}_n \right\rangle\right| \leq 2^{j/2}2^{-2(j + 1)}|H - 1/2| \left(t + \frac{2k + 1 - \theta}{2^{j+1}}\right)^{H - 3/2}.
\end{split}
\end{equation}
Now, as in the proof of Lemma \ref{Sec3Lems1}, we denote by $2^j + \widehat{k_{t, j}}$ the unique $n$ such that $t$ $\in$ $[\frac{\widehat{k_{t, j}}}{2^{j}}, \frac{\widehat{k_{t, j}} + 1}{2^{j}})$. Then, there are two and only two cases:

\textit{Case 1:} \ $\widehat{k_{t, j}}$ $\geq$ $1$. By \eqref{lems2-1-1-02Eq3} we have
\begin{equation*}
\begin{split}
&\left|\left\langle (t + s)^{H - 1/2}, {\cal H}_n \right\rangle\right| \\
&\leq 2^{j/2}2^{-2(j + 1)}|H - 1/2| \left(\frac{\widehat{2k_{t, j}}}{2^{j + 1}} + \frac{2k + 1 - \theta}{2^{j+1}}\right)^{H - 3/2}\\
&\leq 2^{-2jH} \ |H - 1/2| \ 2^{-(H + 1/2)} \ (k + 1)^{H - 3/2}.
\end{split}
\end{equation*}
Thus, there is an absolute constant $D_{2,1}$ $>$ $0$ such that
\begin{equation}
\label{lems2-1-1-02Eq4}
\begin{split}
&\sum_{\{n \text{ at level } j\}}\left\langle (t + s)^{H - 1/2}, {\cal H}_n \right\rangle^2 \leq D_{2, 1}2^{-2jH}
\sum_{\ell = 1}^\infty \left(\frac{1}{\ell}\right)^{3 - 2H}.
\end{split}
\end{equation}

\textit{Case 2:} \ $\widehat{k_{t, j}}$ $=$ $0$. Using \eqref{lems2-1-1-02Eq2} we have
\begin{equation*}
\begin{split}
&\left\langle (t + s)^{H - 1/2}, {\cal H}_{2^j} \right\rangle \\
& = 2^{j/2} \left[\int_{0}^{\frac{1}{2^{j + 1}}} (t + s)^{H - 1/2} ds - \int_{\frac{1}{2^{j + 1}}}^{\frac{2}{2^{j + 1}}} (t + s)^{H - 1/2} ds \right]\\
& = \frac{2^{j/2}}{H + 1/2}\left[\left(\left(t + \frac{1}{2^{j+1}}\right)^{H+1/2} - t^{H+1/2} \right) - \right.\\
&\left.\left(\left(t + \frac{2}{2^{j+1}}\right)^{H+1/2} - \left(t + \frac{1}{2^{j+1}}\right)^{H+1/2}\right)\right]
\end{split}
\end{equation*}
and hence
\begin{equation}
\label{lems2-1-1-02Eq5}
\begin{split}
&\left|\left\langle (t + s)^{H - 1/2}, {\cal H}_{2^j} \right\rangle \right|
\leq \frac{2^{j/2}}{H + 1/2} \left(\frac{2}{2^{j}}\right)^{H+1/2}.
\end{split}
\end{equation}
For $n$ $=$ $2^j$ $+$ $k$ with $k$ $=$ $1$, $\ldots$, $2^j - 1$, by \eqref{lems2-1-1-02Eq3} we have
\begin{equation}
\label{lems2-1-1-02Eq6}
\begin{split}
&\left|\left\langle (t + s)^{H - 1/2}, {\cal H}_n \right\rangle\right| \leq 2^{j/2}2^{-2(j + 1)}|H - 1/2| \left(\frac{2k + 1 - \theta}{2^{j+1}}\right)^{H - 3/2}\\
&< \frac{2^{j/2}}{H + 1/2} 2^{-2(j + 1)}  \left(\frac{k}{2^{j}}\right)^{H -
3/2} \text{ since } |H^2 - 1/4| < 1 \text{ for } H \in (0, 1).
\end{split}
\end{equation}
Putting \eqref{lems2-1-1-02Eq5} and \eqref{lems2-1-1-02Eq6} together, we have
\begin{equation}
\label{lems2-1-1-02Eq7}
\begin{split}
&\sum_{\{n \text{ at level } j\}}\left\langle (t + s)^{H - 1/2}, {\cal H}_n \right\rangle^2 \leq D_{2, 1}2^{-2jH}
\sum_{\ell = 1}^\infty \left(\frac{1}{\ell}\right)^{3 - 2H}.
\end{split}
\end{equation}
Without loss of generality we can let $D_{2, 1}$ be the same absolute constant as in \eqref{lems2-1-1-02Eq4}.

With an argument similar to that for $\left\langle (t + s)^{H - 1/2}, {\cal H}_n \right\rangle$ presented above, we have that there is an absolute constant $D_{2,2}$ $>$ $0$ such that
\begin{equation}
\label{lems2-1-1-02Eq8}
\begin{split}
&\sum_{\{n \text{ at level } j\}}\left\langle s^{H - 1/2}, {\cal H}_n \right\rangle^2 \leq D_{2, 2}2^{-2jH}
\sum_{\ell = 1}^\infty \left(\frac{1}{\ell}\right)^{3 - 2H}.
\end{split}
\end{equation}

The lemma follows from putting  \eqref{lems2-1-1-02Eq1}, \eqref{lems2-1-1-02Eq7} and \eqref{lems2-1-1-02Eq8} together. $\Box$

\begin{lem}
\label{Sec4Lems2}
For any given $H$ $\in$ $(0,1)$ and $q \geq 2$, we have for all $N$ $>$ $1$,
\begin{equation*}
\begin{split}
& \Pr\left\{\sup_{t \in [0,1]  \cap \mathbb Q}\left| I_2(t, H) - W_2(t, H, N) \right| \geq \frac{C_H \sqrt{2 D_2 q}}{\sqrt{H(1 - H)}} \frac{\sqrt{\log N}}{N^H}\right\} \leq \frac{2}{\sqrt{\pi} N^{q}},
\end{split}
\end{equation*}
where $D_2$ is the absolute constant used in Lemma \ref{Sec4Lems1}.
\end{lem}
\textbf{Proof.}
By \eqref{Sec4Eq2} and the consequence of \eqref{Sec4Eq1}, we have
\begin{equation*}
\begin{split}
\Pr \bigcap_{t \in [0, 1] \cap \mathbb Q} & \{(I_2(t, H) - W_2(t, H, N))(\omega) =  \\
& C_H\sum_{n = N + 1}^\infty \langle f_t^{(2)}, \widetilde{\cal H}_n \rangle \int_{-1}^0 \widetilde{\cal H}_n (s) dB_s(\omega) \} = 1.
\end{split}
\end{equation*}
Here $\sum_{n=N+1}^{\infty} \langle f_t^{(2)}, \widetilde{\cal H}_n \rangle\int_{-1}^0 {\cal H}_n(s)dB_s$ is a Gaussian random variable with mean $0$ and variance $\sum_{n=N+1}^{\infty} \langle f_t^{(2)}, \widetilde{\cal H}_n \rangle^2$. The rest of this proof follows the same procedure as in the proof for Lemma \ref{Sec3Lems2}. $\Box$

\section{Approximation of $I_3(t, H)$}
\label{Sec5}
By the time inversion of Bm, we define a Bm $(\widetilde B_s)_{s \in [-1, 0]}$: $\widetilde B_s$ $=$ $s B_{\frac{1}{s}}$ for $s$ $\in$ $[-1, 0)$ and $\widetilde B_0$ $=$ $0$. Consider a family of functions $f_u^{(3)}(v)$ $\in$ $L^2[-1, 0]$ with a parameter $u$ $\in$ $[-1, 0]$: $f_u^{(3)}(v)$ $=$ $1$ if $v$ $\in$ $(u, 0)$; $f_u^{(3)}(v)$ $=$ $0$ otherwise. Let
\begin{equation*}
\begin{split}
g_n(t, H) = \int_{-1}^0  \left((-u^{-1})^{H-3/2} - (t - u^{-1})^{H-3/2}\right) u^{-3} \langle f_u^{(3)}, \widetilde{\cal H}_n \rangle  du.
\end{split}
\end{equation*}
Let $({\cal L}_n^{(3)})_{n \geq 0}$ be the sequence with ${\cal L}_n^{(3)}$ $=$ $\int_{-1}^0 \widetilde {\cal H}_n(s) d\widetilde B_s$. And let ${\cal L}^{*}$ $=$ $B_{-1}$. We define for all $N$ $\geq$ $1$,
\begin{equation}
\label{Sec5Eq1}
\begin{split}
W_3(t, H, N) =  C_H ( (t + 1)^{H-1/2}- 1 ) {\cal L}^{*} - C_H (H - 1/2) \sum_{n = 0}^N g_n(t, H) {\cal L}_n^{(3)}.
\end{split}
\end{equation}

Applying Lemma \ref{Sec3Lems2} to the case where $H$ $=$ $1/2$ and the Haar wavelet $({\cal H}_n)_{n\geq0}$ on $[0,1]$ is replaced by its counterpart  $(\widetilde{\cal H}_n)_{n \geq 0}$ on $[-1, 0]$, we have
\begin{equation}
\label{Sec5Eq2}
\begin{split}
\widetilde B_u = \sum_{n=0}^{\infty} \langle f_u^{(3)}, \widetilde{\cal H}_n \rangle \int_{-1}^0 \widetilde{\cal H}_n(v) d\widetilde B_v \ \text{ almost surely for every } u \in [-1, 0] \cap \mathbb Q.
\end{split}
\end{equation}
A part of Theorem \ref{Thms2} for the case of $H$ $=$ $1/2$ to be proved in next section claims that Lemma \ref{Sec3Lems2} can be extended from discrete time to continuous time. The proof for that part does not involve $I_3(t, H)$ and $I_2(t, H)$ (see Remark on Theorem \ref{Thms2}). We can in this section use the part, i.e., \eqref{Sec5Eq2} can be extended to for every $u$ $\in$ $[-1, 0]$.

\begin{lem}
\label{Sec5Lems1}
There is an absolute constant $D_3$ $>$ $0$ such that for any given $H$ $\in$ $(0,1)$ and $q \geq 2$, we have for all $N$ $>$ $1$,
\begin{equation*}
\begin{split}
\Pr\left\{\sup_{t \in [0,1]  \cap \mathbb Q}\left| I_3(t, H) - W_3(t, H, N) \right| \geq \frac{C_H D_3 \sqrt{q \log N}}{\sqrt{H} N^{1 - H}} \right\} \leq \frac{1}{\sqrt{\pi} N^{q}}.
\end{split}
\end{equation*}
\end{lem}
\textbf{Proof.} Using stochastic integration by parts and the inversion law of Bm, Garz\'on, Gorostiza and Le\'on showed a technical lemma (Lemma 3.1 in \cite{GGL}). By this technical lemma, we have for any fixed $t$ $\in$ $[0, 1]$, almost surely
\begin{equation}
\label{Sec5Lems1Eq1}
\begin{split}
I_3(t, H) & = C_H ( (t + 1)^{H-1/2}- 1 ) B_{- 1} \ - \\
& C_H (H - 1/2) \int_{-1}^0  \left((-u^{-1})^{H-3/2} - (t - u^{-1})^{H-3/2}\right) u^{-3} \widetilde B_u \ du .
\end{split}
\end{equation}
Using the extension of \eqref{Sec5Eq2}, we have for any fixed $t$ $\in$ $[0, 1] \cap \mathbb Q$, almost surely
\begin{equation}
\label{Sec5Lems1Eq2}
\begin{split}
& \int_{-1}^0  \left((-u^{-1})^{H-3/2} - (t - u^{-1})^{H-3/2}\right) u^{-3} \widetilde B_u \ du \\
& = \int_{-1}^0   \sum_{n=0}^{\infty} \left((-u^{-1})^{H-3/2} - (t - u^{-1})^{H-3/2}\right) u^{-3} \langle f_u^{(3)}, \widetilde{\cal H}_n \rangle \int_{-1}^0 \widetilde{\cal H}_n(v) d\widetilde B_v \ du.
\end{split}
\end{equation}

For any fixed $t$ $\in$ $[0, 1] \cap \mathbb Q$, on the right side of \eqref{Sec5Lems1Eq2}, the summation over $n$ and the integration with respect to $du$ are interchangeable. To see this, we regard the summation as a discrete version of integration. By L\'evy's equivalence theorem we have almost surely
\begin{equation}
\label{Sec5Lems1Eq3}
\begin{split}
\sum_{n = 0}^{\infty} \langle f_u^{(3)}, \widetilde{\cal H}_n \rangle^2 \left[ \int_{-1}^0 \widetilde{\cal H}_n(v) d\widetilde B_v \right]^2 = \int_{-1}^0 (f_u^{(3)}(v))^2 dv \left[ \int_{-1}^0 d\widetilde B_v \right]^2 = |u| (\widetilde B_{-1})^2. \\
\end{split}
\end{equation}
And we have for $H$ $\in$ $(0, 1)$, $u$ $\in$ $[-1, 0)$ and $t$ $\in$ $[0, 1]$,
\begin{equation}
\label{Sec5Lems1Eq4}
\begin{split}
& \left| (-u^{-1})^{H-3/2} - (t - u^{-1})^{H-3/2} \right| = |H - 3/2| \int_0^t (s - u^{-1})^{H - 5/2}ds \\
& \leq |H - 3/2| (-u)^{- H + 5/2} \int_0^t ds \leq (3/2) (-u)^{- H + 5/2}.
\end{split}
\end{equation}
By \eqref{Sec5Lems1Eq3} and \eqref{Sec5Lems1Eq4} we have for $H$ $\in$ $(0, 1)$,
\begin{equation*}
\begin{split}
& \int_{-1}^0 \left\{ \sum_{n=0}^{\infty} \left((-u^{-1})^{H-3/2} - (t - u^{-1})^{H-3/2}\right)^2 u^{-6} \langle f_u^{(3)}, \widetilde{\cal H}_n \rangle^2 \left[ \int_{-1}^0 \widetilde{\cal H}_n(v) d\widetilde B_v \right]^2  \right\}^{1/2} du \\
& \leq  \frac{3 | \widetilde B_{-1} | }{2} \int_{-1}^0 (-u)^{-H} du = \frac{3|\widetilde B_{-1}|}{2(1 - H)} < \infty \ \text{ with probability } 1,
\end{split}
\end{equation*}
which implies that the stochastic Fubini theorem is applicable (e.g., see the condition (1.5) in \cite{Ve}). Then it follows from \eqref{Sec5Lems1Eq2} that for any fixed $t$ $\in$ $[0, 1] \cap \mathbb Q$, almost surely
\begin{equation}
\label{Sec5Lems1Eq5}
\begin{split}
& \int_{-1}^0  \left((-u^{-1})^{H-3/2} - (t - u^{-1})^{H-3/2}\right) u^{-3} \widetilde B_u \ du = \sum_{n=0}^{\infty} g_n(t, H) \int_{-1}^0 \widetilde{\cal H}_n(v) d\widetilde B_v .
\end{split}
\end{equation}

Throughout the rest of this proof, we suppose $H$ $\in$ $(0, 1)$ $\setminus$ $\{1/2\}$. Consider a family of functions $f_x^{(4)}(s)$ $\in$ $L^2[0, 1]$ with a parameter $x$ $\in$ $[0, 1]$: $f_x^{(4)}(s)$ $=$ $1$ if $s$ $\in$ $(0, x)$; $f_x^{(4)}(s)$ $=$ $0$ otherwise. Replacing $x$ by $-u$ and $(\widetilde{\cal H}_n)_{n \geq 0}$ by $({\cal H}_n)_{n \geq 0}$, we have
\begin{equation}
\label{Sec5Lems1Eq6}
\begin{split}
g_n(t, H) & = \int_0^1 \left((t + x^{-1})^{H-3/2} - (x^{-1})^{H-3/2} \right) x^{-3} \langle f_x^{(4)}, {\cal H}_n \rangle dx.
\end{split}
\end{equation}
Recall the two conventions: $n \in \mathbb Z^+$ is said to be at level $j$ if $n$ $=$ $2^j$ $+$ $k$ with $j$ $\geq$ $0$ and $0$ $\leq$ $k$ $<$ $2^j$; and the interval $[\frac{k}{2^j}, \frac{k+1}{2^j})$ is meant to be $[\frac{k}{2^j}, \frac{k+1}{2^j}]$ when $\frac{k+1}{2^j}$ $=$ $1$. For $n$ $=$ $2^j$ $+$ $k$, let $g_{j, k}(t, H)$ $=$ $\int_{\frac{k}{2^j}}^{\frac{k+1}{2^j}} \left((t + x^{-1})^{H-3/2} - (x^{-1})^{H-3/2} \right) x^{-3} \langle f_x^{(4)}, {\cal H}_n \rangle dx$. For notional simplicity, let $G_{t, H}(x)$ $=$ $\left((t + x^{-1})^{H-3/2} - (x^{-1})^{H-3/2} \right) x^{-3}$. We have
\begin{equation}
\label{Sec5Lems1Eq7}
\begin{split}
g_{j, k}(t, H) = & \ \int_{\frac{2k}{2^{j + 1}}}^{\frac{2k + 1}{2^{j + 1}}} G_{t, H}(x) \int_0^x {\cal H}_n(y) dy \, dx + \int_{\frac{2k + 1}{2^{j + 1}}}^{\frac{2k + 2}{2^{j + 1}}} G_{t, H}(x) \int_0^x {\cal H}_n(y) dy \, dx \\
= & \ \  2^{j/2} \int_{\frac{2k}{2^{j + 1}}}^{\frac{2k + 1}{2^{j + 1}}} G_{t, H}(x) x \, dx - 2^{j/2} \int_{\frac{2k + 1}{2^{j + 1}}}^{\frac{2k + 2}{2^{j + 1}}} G_{t, H}(x) x \, dx \\
& - 2^{j/2} \int_{\frac{2k}{2^{j + 1}}}^{\frac{2k + 1}{2^{j + 1}}} G_{t, H}(x) \frac{2k}{2^{j + 1}} \, dx + 2^{j/2} \int_{\frac{2k + 1}{2^{j + 1}}}^{\frac{2k + 2}{2^{j + 1}}} G_{t, H}(x) \frac{2k + 2}{2^{j + 1}} \, dx.
\end{split}
\end{equation}

For the first two terms on the rightmost side of \eqref{Sec5Lems1Eq7}, we have $\int_a^b G_{t, H}(x) x \, dx$ $=$ $\frac{1}{H - 1/2}\left(y^{H - 1/2}\right.$ $-$ $\left.\left.(t + y)^{H - 1/2} \right)\right|_{y = 1/a}^{y = 1/b}$ for $b$, $a$ $>$ $0$. We denote $2^{j/2} \int_{\frac{2k}{2^{j + 1}}}^{\frac{2k + 1}{2^{j + 1}}} G_{t, H}(x) x \, dx$ $-$ $2^{j/2} \int_{\frac{2k + 1}{2^{j + 1}}}^{\frac{2k + 2}{2^{j + 1}}} G_{t, H}(x) x \, dx$ by $\widetilde h_{t, H, j, \, k}$. In the case of $k$ $=$ $0$, we have
\begin{equation*}
\begin{split}
& \widetilde h_{t, H, j, \, 0} = \frac{2^{j/2}}{ H - 1/2 } \left[ \left. \left( y^{H - 1/2} - (t + y)^{H - 1/2} \right) \right|_{y = \infty}^{y = 2^{j + 1}} - \left( y^{H - 1/2} - \right. \right. \\
& \left. \left. \left. (t + y)^{H - 1/2} \right) \right|_{y = 2^{j + 1}}^{y = 2^j} \right]  = \frac{2^{j/2}}{ H - 1/2 } \left[  2^{(j + 1)(H - 1/2) + 1} \left(1  - (\frac{t}{2^{j + 1}} + 1)^{H - 1/2} \right) \right. \\
& - \left. 2^{j(H - 1/2)} \left(1  - (\frac{t}{2^{j}} + 1)^{H - 1/2} \right) \right], \\
\end{split}
\end{equation*}
which implies for $t$ $\in$ $[0. 1]$,
\begin{equation}
\label{Sec5Lems1Eq8}
\begin{split}
|\widetilde h_{t, H, j, \, 0}| \leq \frac{D_{3, 1}^{*}}{ 2^{j(1 - H)}} \ \text{ with an absolute constant } D_{3, 1}^{*} > 0.
\end{split}
\end{equation}
In the case of $k$ $>$ $0$, we have
\begin{equation}
\label{Sec5Lems1Eq9}
\begin{split}
& \widetilde h_{t, H, j, \, k} = \frac{2^{j/2}}{ H - 1/2 } \left[ \left. \left( y^{H - 1/2} - (t + y)^{H - 1/2} \right) \right|_{y = \frac{2^{j + 1}}{2k}}^{y = \frac{2^{j + 1}}{2k + 1}} - \left( y^{H - 1/2} - \right. \right. \\
& \left. \left. \left. (t + y)^{H - 1/2} \right) \right|_{y = \frac{2^{j + 1}}{2k + 1}}^{y = \frac{2^{j + 1}}{2k + 2}} \right] =  \frac{ 2^{jH} \, 2^{H - 1/2} }{ H - 1/2 } \left[ \left. 2 \left( 1 - (\frac{t(2k + 1)}{2^{j + 1}} + 1)^{H - 1/2} \right) \right. \right. \\
& - \left. \left( 1 - (\frac{t2k}{2^{j + 1}} + 1)^{H - 1/2} \right) - \left( 1 - ( \frac{t(2k + 2)}{2^{j + 1}} + 1)^{H - 1/2} \right) \right].
\end{split}
\end{equation}
For the rightmost side of \eqref{Sec5Lems1Eq9} we introduce a function $\widetilde w$ of $h$: $\widetilde w(h)$ $=$ $2 \, \widetilde g(x_0)$ $-$ $\widetilde g(x_0 + h)$ $-$ $\widetilde g(x_0 - h)$ where $\widetilde g( x )$ $=$ $1 - (1 + \frac{t(2k + 1 + x)}{2^{j + 1}})^{H - 1/2}$ and $x_0$ $=$ $0$. Then we have $\widetilde h_{t, H, j, \, k}$ $=$ $\frac{ 2^{jH} \, 2^{H - 1/2} }{ H - 1/2 } \widetilde w(1)$. By Taylor's expansion we have
\begin{equation*}
\begin{split}
&\widetilde w(h) = \widetilde w(0) + \frac{\widetilde w'(0)}{1!}h + \frac{\widetilde w''(\theta h)}{2!}h^2 \ \text{ (for some } 0 < \theta < 1)\\
& = \frac{\widetilde w''(\theta h)}{2!}h^2 \ \text{(since } \widetilde w(0) = \widetilde w'(0) = 0), \ \text{ where } \widetilde w''(x) = \frac{(H - 1/2) (H - 3/2) t^2}{2^{2(j + 1)}} \times \\
& \left[ (1 + \frac{t(2k + 1 + x)}{2^{j + 1}})^{H - 5/2} +  (1 + \frac{t(2k + 1 - x)}{2^{j + 1}})^{H - 5/2} \right].
\end{split}
\end{equation*}
Hence, we have an absolute constant $D_{3, 2}^{*}$ $>$ $0$ so that for $n$ $=$ $2^j$ $+$ $k$ with $0$ $<$ $k$ $< 2^j$,
\begin{equation}
\label{Sec5Lems1Eq10}
\begin{split}
& | \widetilde h_{t, H, j, \, k} | = \left| \frac{ 2^{jH} \, 2^{H - 1/2} }{ H - 1/2 } \widetilde w(1) \right| = \left| \frac{ 2^{jH} \, 2^{H - 1/2} }{ H - 1/2 } \frac{\widetilde w''( \theta )}{2!} \right| \leq \frac{D^{*}_{3, 2}}{2^{j(2 - H)}| H - 1/2 |}  .
\end{split}
\end{equation}

Using the same method, we estimate the last two terms on the rightmost side of \eqref{Sec5Lems1Eq7}. Denote $\left( - 2^{j/2} \int_{\frac{2k}{2^{j + 1}}}^{\frac{2k + 1}{2^{j + 1}}} G_{t, H}(x) \frac{2k}{2^{j + 1}} \, dx \right.$ $+$ $\left.  2^{j/2} \int_{\frac{2k + 1}{2^{j + 1}}}^{\frac{2k + 2}{2^{j + 1}}} G_{t, H}(x) \frac{2k + 2}{2^{j + 1}} \, dx \right)$ by $\widehat h_{t, H, j, \, k}$. In the case of $k$ $=$ $0$, we have $\widehat h_{t, H, j, \, 0}$ $=$ $2^{- j/2} \int_{\frac{1}{2^{j + 1}}}^{\frac{2}{2^{j + 1}}} G_{t, H}(x) \, dx$. Then, using \eqref{Sec5Lems1Eq4} we have an absolute constant $D^{*}_{3, 3}$ $>$ $0$ such that
\begin{equation}
\label{Sec5Lems1Eq11}
\begin{split}
& | \widehat h_{t, H, j, \, 0} | \leq \frac{3 \times 2^{- j/2}}{2} \int_{\frac{1}{2^{j + 1}}}^{\frac{2}{2^{j + 1}}} x^{- H - 1/2} dx \leq  \frac{D^{*}_{3, 3}}{2^{j(1 - H)}}.
\end{split}
\end{equation}
For the case of $k$ $>$ $0$, we have for $b$ $>$ $a$ $>$ $0$,
\begin{equation}
\label{Sec5Lems1Eq12}
\begin{split}
& \int_a^b G_{t, H}(x) dx = \int_{\frac{1}{b}}^{\frac{1}{a}} \left( (t + u)^{H - 3/2} - u^{H - 3/2} \right) u \ du  \\
= & \ \left( \frac{1}{a} \right)^{H + 1/2} \left[ \frac{(at + 1)^{H - 1/2}}{H - 1/2} - \frac{1}{H + 1/2} - \frac{(at + 1)^{H + 1/2}}{(H - 1/2)(H + 1/2)} \right]  \\
& - \left( \frac{1}{b} \right)^{H + 1/2} \left[ \frac{(bt + 1)^{H - 1/2}}{H - 1/2} - \frac{1}{H + 1/2} - \frac{(bt + 1)^{H + 1/2}}{(H - 1/2)(H + 1/2)} \right] .
\end{split}
\end{equation}
We introduce a function: $\widehat w(h)$ $=$ $\widehat g(x_0 + h)$ $+$ $\widehat g(x_0 - h)$ $-$ $2 \widehat g(x_0)$ where $x_0$ $=$ $0$ and $\widehat g(x)$ $=$ $\frac{2k + 1 + x}{2} \left( \frac{2^{j + 1}}{2k + 1 + x} \right)^{H + 1/2}
\left[  \frac{(\frac{(2k + 1 + x)t}{2^{j + 1}} + 1)^{H - 1/2}}{H - 1/2} - \frac{1}{H + 1/2} - \frac{(\frac{(2k + 1 + x)t}{2^{j + 1}} + 1)^{H + 1/2}}{(H - 1/2)(H + 1/2)} \right]$. We denote by $f(x)$ the third factor $[\cdot]$ in $\widehat g(x)$. By Taylor's expansion we have for some $0$ $<$ $\theta$ $<$ $1$,
\begin{equation}
\label{Sec5Lems1Eq13}
\begin{split}
&\widehat w(h) = \widehat w(0) + \frac{\widehat w'(0)}{1!}h + \frac{\widehat w''(\theta h)}{2!}h^2 = \frac{\widehat w''(\theta h)}{2!}h^2, \ \ \ \widehat w''(x) = 2^{j(H + 1/2) + H - 1/2} \times \\
&  \left\{ \left[ \frac{(H - 1/2)(H + 1/2) f(x)}{ (2k + 1 + x)^{H + 3/2} } - \frac{(2 H - 1)f'(x)}{(2k + 1 + x)^{H + 1/2}} + \frac{f''(x)}{(2k + 1 + x)^{H - 1/2}}  \right] + \right. \\
& \left. \left[ \frac{(H - 1/2)(H + 1/2) f(-x)}{ (2k + 1 - x)^{H + 3/2} }  - \frac{(2 H - 1)f'(- x)}{(2k + 1 - x)^{H + 1/2}} + \frac{f''(- x)}{(2k + 1 - x)^{H - 1/2}}  \right] \right\}.
\end{split}
\end{equation}
By \eqref{Sec5Lems1Eq12} and \eqref{Sec5Lems1Eq13} we have $\widehat h_{t, H, j, \, k}$ $=$ $2^{-j/2} \widehat w(1)$ $=$ $\frac{2^{-j/2} \widehat w''( \theta )}{2!}$. Then, using calculus we have an estimate as follows (for which details are in an appendix available upon a request sent to the corresponding author): there is an absolute constant $D_{3, 4}^{*}$ $>$ $0$ such that for $n$ $=$ $2^j$ $+$ $k$ with $0$ $<$ $k$ $<$ $2^j$,
\begin{equation}
\label{Sec5Lems1Eq14}
\begin{split}
| \widehat h_{t, H, j, \, k} | & = \left| \frac{2^{-j/2} \widehat w''( \theta )}{2!} \right| \leq  \frac{D_{3, 4}^{*}}{2^{j(1- H)}} \left( \frac{1}{k + 1} \right)^{H + 1/2}.
\end{split}
\end{equation}

Now, putting \eqref{Sec5Lems1Eq6},  \eqref{Sec5Lems1Eq7}, \eqref{Sec5Lems1Eq8}, \eqref{Sec5Lems1Eq10}, \eqref{Sec5Lems1Eq11}, and \eqref{Sec5Lems1Eq14} together, we have an estimate as follows: there is an absolute constant $D_{3, 1}$ $>$ $0$ such that
\begin{equation}
\label{Sec5Lems1Eq15}
\begin{split}
& \sum_{\{ n \text{ at level } j \}} [g_n(t, H)]^2 \leq \frac{D_{3, 1}}{ 2^{2j(1 - H)} (H - 1/2)^2 } \sum_{\ell = 1}^{\infty} \left( \frac{1}{\ell} \right)^{2H + 1}.
\end{split}
\end{equation}
Then, by \eqref{Sec5Lems1Eq1}, \eqref{Sec5Lems1Eq2}, \eqref{Sec5Lems1Eq5}, and \eqref{Sec5Lems1Eq15} we can use a procedure similar to the proof of Lemma \ref{Sec3Lems1} and then use a procedure similar to the proof for Lemma \ref{Sec3Lems2} to complete a proof for this lemma. $\Box$

\noindent
\textbf{Remark.} In the above proof, the time inversion of Bm adds a factor $u^{-1}$ to the integrand $\left((-u^{-1})^{H-3/2} - (t - u^{-1})^{H-3/2}\right) u^{-3}$ in the second term on the right side of \eqref{Sec5Lems1Eq1} where $u^{-2}$ in $u^{-3}$ is from change of variable. Denote the integrand by $\cal Q$. We have $\cal Q$ $\sim$ $u^{- H - 1/2}$ as $u$ $\rightarrow$ $0$. The exponent in $u^{- H - 1/2}$ causes the convergence rate $O(N^{- (1 - H)} \sqrt {\log N})$ of $W_3(t, H)$ to $I_3(t, H)$. For $H$ $\in$ $(0, 1/2)$, it is faster than the convergence rate $O(N^{- H} \sqrt {\log N})$ of $W_1(t, H)$ to $I_1(t, H)$ as well as $W_2(t, H)$ to $I_2(t, H)$. But, for $H$ $\in$ $(1/2, 1)$, the convergence rate caused by the exponent becomes slow, which reflects an impact of the long-range dependence of an fBm of $H$ $\in$ $(1/2, 1)$.

\section{Approximation of fBm}
\label{Sec6}
In $(\Omega, {\cal F}, \Pr)$ we define for $t$ $\in$ $[0, 1] \cap \mathbb Q$ and $q$ $\geq$ $2$,
\begin{equation*}
\begin{split}
W(t, H, N) = & \ W_1(t, H, N) + W_2(t, H, N) + W_3(t, H, N).
\end{split}
\end{equation*}
By Lemma \ref{Sec3Lems2}, Lemma \ref{Sec4Lems2}, Lemma \ref{Sec5Lems1} and $H$ $>$ $1$ $-$ $H$ for $H$ $\in$ $(1/2, 1)$, we have
\begin{thm}
\label{Thms1}
There are absolute constants $C_{1, 1}$, $C_{1, 2}$, $C_{2, 1}$, $C_{2, 2}$ $>$ $0$ such that for any given $H$ $\in$ $(0,1/2]$ and $q \geq 2$, we have for all $N$ $>$ $1$,
\begin{equation}
\label{Thms1Eq1}
\begin{split}
&\Pr\left\{\sup_{t\in[0,1] \cap \mathbb Q}\left|B_t^{(H)} - W(t, H, N) \right| \geq \frac{C_{1, 1} \sqrt q}{\sqrt{H(1 - H)}} \frac{\sqrt{\log N}}{N^H}\right\} \leq \frac{C_{1, 2}}{N^q},
\end{split}
\end{equation}
and for any given $H$ $\in$ $(1/2, 1)$ and $q \geq 2$, we have for all $N$ $>$ $1$,
\begin{equation}
\label{Thms1Eq2}
\begin{split}
&\Pr\left\{\sup_{t\in[0,1] \cap \mathbb Q}\left|B_t^{(H)} - W(t, H, N) \right| \geq \frac{C_{2, 1} \sqrt q}{\sqrt{H(1 - H)}} \frac{\sqrt{\log N}}{N^{1 - H}}\right\} \leq \frac{C_{2, 2}}{N^q}. \ \Box
\end{split}
\end{equation}
\end{thm}

With respect to a H\"older continuous version of an fBm, Theorem \ref{Thms1} can be extended from discrete time $t$ $\in$ $[0, 1]$ $\cap$ $\mathbb Q$ to continuous time $t$ $\in$ $[0, 1]$.
\begin{thm}
\label{Thms2}
An fBm $(B_t^{(H)})_{t \in [0, 1]}$ of $H$ $\in$ $(0, 1)$ has a wavelet-based almost sure uniform expansion as follows: In $(\Omega, {\cal F}, \Pr)$ we have for $t$ $\in$ $[0, 1]$, with probability $1$
\begin{equation*}
\begin{split}
B_t^{(H)} = & \ C_H \sum_{n = 0}^{\infty} \langle f_t^{(1)}, {\cal H}_n  \rangle {\cal L}^{(1)}_n +  C_H  \sum_{n = 0}^{\infty} \langle f_t^{(2)}, \widetilde{\cal H}_n  \rangle {\cal L}^{(2)}_n + C_H ( (t + 1)^{H-1/2}- 1 ) {\cal L}^* \\
& + C_H (H - 1/2) \sum_{n = 1}^{\infty} g_n(t, H) {\cal L}^{(3)}_n, \ \text{ where}
\end{split}
\end{equation*}
$\langle f_t^{(1)}, {\cal H}_n  \rangle$ and $({\cal L}^{(1)}_n)_{n \geq 0}$, $\langle f_t^{(2)}, \widetilde{\cal H}_n  \rangle$ and $({\cal L}^{(2)}_n)_{n \geq 0}$, and ${\cal L}^*$, $g_n(t, H)$ and $({\cal L}^{(3)}_n)_{n \geq 0}$ are the same as in \eqref{Sec3Eq2}, \eqref{Sec4Eq2}, and \eqref{Sec5Eq1}, respectively. Regarding $\sum_{n = 0}^{\infty}$ as $\lim_{N \rightarrow \infty} \sum_{n = 0}^{N}$, the convergence rates are $O(N^{-H} \sqrt{ \log N})$ for $H$ $\in$ $(0, 1/2]$ and $O(N^{-(1 - H)} \sqrt{ \log N})$ for $H$ $\in$ $(1/2, 1)$, as expressed by \eqref{Thms1Eq1} and \eqref{Thms1Eq2} respectively.
\end{thm}
Recall that by \eqref{Sec2Eq1} we write $B_t^{(H)}$ as $I_1(t, H)$ $+$ $I_2(t, H)$ $+$ $I_3(t, H)$ which then are approximated by $W_1(t, H)$, $W_2(t, H)$, and $W_3(t, H)$ separately. We below provide a proof for the extension of the approximation of $I_1(t, H)$ by $W_1(t, H)$ from $t$ $\in$ $[0, 1] \cap \mathbb Q$  to $t$ $\in$ $[0, 1]$ in the case of $H$ $\in$ $(0, 1/2]$. Proofs for $H$ $\in$ $(1/2, 1)$ and all other cases, including the extension of the approximation of $I_2(t, H)$ by $W_2(t, H)$ as well as $I_3(t, H)$ by $W_3(t, H)$, can be carried out in a similar way.

\noindent
\textbf{Proof for Theorem \ref{Thms2}.} Using $\alpha$ $\in$ $\mathbb Z^+$ as a parameter, we let $[0, 1]$ $=$ $\bigcup_{i = 1}^{16^{\alpha}} [\frac{i - 1}{16^{\alpha}}, \frac{i}{16^{\alpha}}) \cup \{1\}$. We denote $\left( C_H \sum_{n=0}^N \langle f_{\frac{i - 1}{16^{\alpha}}}^{(1)}, {\cal H}_n \rangle \int_0^1 {\cal H}_n(s) dB_s \right)$ by ${\cal M}_1 (N, \alpha, i)$. For $t^*$ $\in$ $[\frac{i - 1}{16^{\alpha}}, \frac{i}{16^{\alpha}})$ $\setminus$ $\mathbb Q$, we denote $\left( C_H \sum_{n=0}^N \langle f_{t*}^{(1)}, {\cal H}_n \rangle \int_0^1 {\cal H}_n(s) dB_s \right.$ $-$ $\left. C_H \sum_{n=0}^N \langle f_{\frac{i - 1}{16^{\alpha}}}^{(1)}, {\cal H}_n \rangle \int_0^1 {\cal H}_n(s) dB_s \right)$ by $Q_1(t^*, N, \alpha, i)$. By Lemma \ref{Sec3Lems2} we have for $q \geq 2$ and for all $N$ $>$ $1$,
\begin{equation}
\label{Thms2Eq1}
\begin{split}
\Pr \left\{ \sup_{\alpha \in \mathbb Z^+, \, 1 \leq i \leq 16^{\alpha}} \left|{\cal M}_1(N, \alpha, i) - B_{\frac{i - 1}{16^{\alpha}}}^{(H)} \right| \geq \frac{C_H\sqrt{2 D_1 q}}{\sqrt{H(1 - H)}} \frac{\sqrt{\log N}}{N^H}\right\} \leq \frac{1}{\sqrt{\pi} N^q}.
\end{split}
\end{equation}
Recall the H\"older continuity of fBm. Almost surely, a sample path $B_t^{(H)}(\omega)$ ($t$ $\in$ $[0, 1]$) is H\"older continuous of order $\beta H$ for $\beta$ $\in$ $(0, 1)$ where $\beta$ cannot be $1$ by the law of iterated logarithm \cite{Ar}. We choose $\beta$ close to $1$, having
\begin{equation}
\label{Thms2Eq2}
\begin{split}
\Pr \left\{ \sup_{\alpha \in \mathbb Z^+, \, 1 \leq i \leq 16^{\alpha}} \left\{ \left|B_{\frac{i - 1}{16^{\alpha}}}^{(H)} - B_{t^*}^{(H)}\right| : t^* \in [\frac{i - 1}{16^{\alpha}}, \frac{i}{16^{\alpha}}) \setminus \mathbb Q \right\} \leq \frac{M}{16^{\alpha \beta}} \right\} = 1
\end{split}
\end{equation}
where $M$ $>$ $0$ is a constant depending only on the chosen $\beta$.

$Q_1(t^*, N, \alpha, i)$ is a Gaussian random variable with mean $0$ and variance $\left( C_H^2 \sum_{n = 0}^N \right.$ $(\int_0^1 (f_{t*}^{(1)}(s)$ $-$ $\left. f_{\frac{i - 1}{16^{\alpha}}}^{(1)}(s)){\cal H}_n(s) ds)^2 \right)$. We estimate the variance. Without loss of generality we suppose $\alpha$ $>$ $\log_ 2 N$. Then, for all $1$ $\leq$ $n$ $\leq$ $N$ and $0$ $\leq$ $k$ $<$ $2^j$, one and only one of the following three cases occurs: $[\frac{i - 1}{16^{\alpha}}, \frac{i}{16^{\alpha}})$ $\subset$ $[\frac{2k}{2^{j + 1}}, \frac{2k + 1}{2^{j + 1}})$; $[\frac{i - 1}{16^{\alpha}}, \frac{i}{16^{\alpha}})$ $\subset$ $[\frac{2k + 1}{2^{j + 1}}, \frac{2k + 2}{2^{j + 1}})$; or $[\frac{i - 1}{16^{\alpha}}, \frac{i}{16^{\alpha}})$ $\cap$ $[\frac{k}{2^j}, \frac{k + 1}{2^j})$ $=$ $\emptyset$. Then by calculus we have an estimate (for which details are in an appendix available upon a request sent to the corresponding author): for any $t^*$ $\in$ $[\frac{i - 1}{16^{\alpha}}, \frac{i}{16^{\alpha}})$ $\setminus \mathbb Q$,
\begin{equation}
\label{Thms2Eq4}
\begin{split}
& \Pr \left\{ | Q_1(t^*, N, \alpha, i) | > \frac{\sqrt{2 G_1}}{(\sqrt 2)^{\alpha}} \right\} \leq \frac{1}{\sqrt \pi} e^{- 2^{\alpha}} \ \text{with an absolute constant } G_1 > 0.
\end{split}
\end{equation}

Consider the H\"older continuous version described in \eqref{Thms2Eq2} over every time interval $t$ $\in$ $[\frac{i - 1}{16^{\alpha}}, \frac{i}{16^{\alpha}})$. Then, by \eqref{Thms2Eq4} we have an absolute constant $G$ $>$ $0$ such that
\begin{equation}
\label{Thms2Eq5}
\begin{split}
\Pr \bigcup_{i = 1}^{16^{\alpha}} & \left\{ \sup_{t^* \in [\frac{i - 1}{16^{\alpha}}, \frac{i}{16^{\alpha}}) \setminus \mathbb Q} \left| {\cal M}^*(t^*, N, \alpha, i) -  B_{t^*}^{(H)} \right| > \frac{G}{(\sqrt 2)^{\alpha}} \ + \right.\\
& \left. \frac{C_H\sqrt{2 D_1 q}}{\sqrt{H(1 - H)}} \frac{\sqrt{\log N}}{N^H} + \frac{M}{16^{\alpha \beta}} \right\} \leq \frac{3 \times 16^{\alpha} }{\sqrt \pi} e^{- 2^{\alpha}} + \frac{1}{\sqrt{\pi} N^q}.
\end{split}
\end{equation}
\eqref{Thms2Eq5} holds for all $\alpha$ $>$ $\log_2 N$. Given $H$ $\in$ $(0,1/2]$ and $q \geq 2$, by \eqref{Thms2Eq5} and \eqref{Thms1Eq1} we have for all $N$ $>$ $1$,
\begin{equation*}
\begin{split}
& \Pr \left\{ \sup_{t \in [0, 1]} \left| B_t^{(H)} - C_H \sum_{n = 0}^{N} \langle f_t^{(1)}, {\cal H}_n  \rangle {\cal L}_n^{(1)}  \right| \geq \frac{\sqrt q (C_H\sqrt{2 D_1} + C_{1, 1})}{\sqrt{H(1 - H)}} \frac{\sqrt{\log N}}{N^H} \right\} \\
& \leq \frac{C_{1, 2} + 1}{\sqrt{\pi} N^q}. \ \ \ \Box
\end{split}
\end{equation*}

\noindent
\textbf{Remark.} The above proof shows that $I_1(t, 1/2)$, which is a Bm, has an almost sure and uniform expansion $\sum_{n = 0}^{\infty} \langle f_t^{(1)}, {\cal H}_n  \rangle \int_0^1 {\cal H}_n (s) dB_s$ for $t$ $\in$ $[0, 1]$. The proof does not involve $I_2$ and $I_3$, which justifies our use of this expansion in the previous section.

\section{A parallel algorithm for the approximation}
\label{Sec7}
We give a mathematical description of an algorithm to demonstrate how a sample path of fBm can be generated in parallel over time. The reader who is in interested in parallel algorithms is referred to \cite{Le}. Theorem \ref{Thms2} implies that a sample path $B_t^{(H)}(\omega) :$ $t \in [0, 1]$ $\mapsto$ $\mathbb R$ can almost surely and uniformly be approximated by
\begin{equation}
\label{Sec7Eq1}
\begin{split}
B_t^{(H)} (\omega) \approx & \ C_H \sum_{n = 0}^N \langle f_t^{(1)}, {\cal H}_n  \rangle {\cal L}^{(1)}_n (\omega) +  C_H  \sum_{n = 0}^N \langle f_t^{(2)}, \widetilde{\cal H}_n  \rangle {\cal L}^{(2)}_n  (\omega) \\
& + C_H ( (t + 1)^{H-1/2}- 1 ) {\cal L}^* (\omega) + C_H (H - 1/2) \sum_{n = 1}^N g_n(t, H) {\cal L}^{(3)}_n (\omega).
\end{split}
\end{equation}
Hence, given any time instances $t_1$, $\ldots$, $t_{\ell}$ $\in$ $[0, 1]$, we can compute approximations of $B_{t_1}^{(H)}(\omega)$, $\ldots$, $B_{t_{\ell}}^{(H)}(\omega)$ as follows. Make $(3N + 4)$ independent observations of a normal distribution ${\cal N}(0, 1)$. Denote the results from the first $(N + 1)$ observations by ${\cal L}_n^{(1)}(\omega)$, $n$ $=$ $0$, $\ldots$, $N$; denote the results from the second $(N + 1)$ observations by ${\cal L}_n^{(2)}(\omega)$, $n$ $=$ $0$, $\ldots$, $N$; denote the results from the third $(N + 1)$ observations by ${\cal L}_n^{(3)}(\omega)$, $n$ $=$ $0$, $\ldots$, $N$; and denote the result from the last observation by ${\cal L}_n^{*}(\omega)$. Then, using \eqref{Sec7Eq1} we compute approximations of $B_{t_1}^{(H)}(\omega)$, $\ldots$, $B_{t_{\ell}}^{(H)}(\omega)$ separately in an arbitrarily chosen order of $t_1$, $\ldots$, $t_{\ell}$. It means that the $\ell$ approximations can be carried out in parallel over time $t_1$, $\ldots$, $t_{\ell}$ $\in$ $[0, 1]$ on multiple (e.g., $\ell$ in the ideal case) processors available in today's computer systems.

By \eqref{Sec7Eq1} we can see that the number $\ell$ of time instances is not related to $N$, the number of approximation step. Given $N$, we can decide at what time instances $t_1$, $\ldots$, $t_{\ell}$ $\in$ $[0, 1]$ we want to find approximations of $B_{t_1}^{(H)}(\omega)$, $\ldots$, $B_{t_{\ell}}^{(H)}(\omega)$. And the accuracy of such approximations is determined by $N$, as shown by the deviation bounds \eqref{Thms1Eq1} and \eqref{Thms1Eq2} respectively for cases of $H$ $\in$ $(0, 1/2]$ and $H$ $\in$ $(1/2, 1)$. Given time instances $t_1$, $\ldots$, $t_{\ell}$ $\in$ $[0, 1]$, we can decide the number $N$ of approximation steps to ensure the accuracy of approximation by the above two deviation bounds.

By using the Mandelbrot - van Ness representation and Haar wavelets, the coefficients on the right side of \eqref{Sec7Eq1}, i.e., $\langle f_t^{(1)}, {\cal H}_n  \rangle$, $\langle f_t^{(2)}, {\cal H}_n  \rangle$ and $g_n(t, H)$ are easy to compute.

\ack
The authors are grateful to the referee for instructive comments which led to significant improvements of the original version of this paper. The authors thank the editor for granting them opportunities to improve the paper.

%
%
%
%

\end{document}